\documentclass[12pt]{article}

\newtheorem{definition}{Definition}

\newcommand{\Z}{{\mathbb Z}}

\newcommand{\JKLM}[4]{
\begin{array}{lll}
J & = & \{ #1 \} \\
K & = & \{ #2 \} \\
L & = & \{ #3 \} \\
M & = & \{ #4 \} \\
\end{array}
}

\usepackage{epsfig}
\usepackage{url}
\usepackage{amssymb}
\usepackage{amsmath}
\usepackage{amsfonts}

\def\Dbar{\leavevmode\lower.6ex\hbox to 0pt{\hskip-.23ex \accent"16\hss}D}

\begin{document}

{\bf\LARGE
\begin{center}
Some new orders of Hadamard and Skew-Hadamard matrices
\end{center}
}

{\Large
\begin{center}
Dragomir {\v{Z}}. {\Dbar}okovi{\'c}\footnote{University of Waterloo,
Department of Pure Mathematics, Waterloo, Ontario, N2L 3G1, Canada
e-mail: \url{djokovic@math.uwaterloo.ca}},
Oleg Golubitsky\footnote{Google, Inc. 151 Charles West, Kitchener, Ontario, N2G 1H6, Canada, e-mail: \url{oleg.golubitsky@gmail.com}},
Ilias S.
Kotsireas\footnote{Wilfrid Laurier University, Department of Physics
\& Computer Science, Waterloo, Ontario, N2L 3C5, Canada, e-mail:
\url{ikotsire@wlu.ca}}
\end{center}
}

\begin{abstract}
\noindent
We construct Hadamard matrices of orders $4 \cdot 251 = 1004$ and $4 \cdot 631 = 2524$, and skew-Hadamard matrices of orders $4\cdot 213 = 852$ and $4 \cdot 631 = 2524$. As far as we know,
such matrices have not been constructed previously.
The constructions use the Goethals-Seidel array, suitable supplementary difference sets on a cyclic group
and a new efficient matching algorithm based on hashing techniques.
\end{abstract}

\section{Introduction}

\noindent There are only $13$ integers $v<500$ for which no Hadamard matrices of order $4v$ are known, see
\cite{Colbourn:Dinitz:Handbook:2006,Djokovic:C:2008_HM,
Horadam:2012}. In this paper we remove the integer $v = 251$ from this list by constructing a Hadamard matrix of order $4 \cdot 251 =1004$
by constructing cyclic supplementary difference sets (SDS) with parameters $(251;125,120,115,115;224)$.
Using the same method we also construct an SDS with parameters
$(631; 315; 330; 330; 330; 647)$ leading to Hadamard matrix of order $4 \cdot 631$ of skew-type.
Our computation also found a suitable SDS which leads
to a Hadamard matrix of skew-type of order $4 \cdot 213$.



We recall the definition of SDS. Let $\Z_v=\Z/v\Z$ denote the ring of integers modulo $v$, and let $|X|$ denote the
cardinality of a finite set $X$. Let $k_1,\ldots,k_t$ and $\lambda$ be positive integers such that
$\lambda (v-1)=\sum k_i(k_i-1)$.

\begin{definition}
We say that $X_1,\ldots,X_t\subseteq\Z_v$ are {\em supplementary difference sets}  with parameters
$(v;k_1,\ldots,k_t;\lambda)$, if $|X_i|=k_i$ for each index $i$ and for every nonzero element $c\in\Z_v$
there are exactly $\lambda$ ordered pairs $(a,b)$ such that $a-b=c \pmod{v}$ and
$\{a,b\}\subseteq X_i$ for some $i\in\{1,2,\ldots,t\}$.
\end{definition}

\noindent The existence of an SDS with parameters $(v;k_1,k_2,k_3,k_4;\lambda)$ with $\lambda=\sum k_i-v$ implies the existence of a Hadamard matrix of order $4v$. If moreover $v$ is odd and $X_1$ is a skew subset of $\Z_v$, i.e., for each nonzero $i\in\Z_v$ exactly one of the integers $i$ and $v-i$ belongs to $X_1$, then one can use such SDS to construct a skew-Hadamard matrix of order $4v$, see \cite{Seberry:Yamada:1992}. \\

\noindent All SDSs in this paper were constructed by using the well-known method of taking base blocks to be unions of orbits of an automorphism group of the underlying cyclic group \cite{Gysin:Seberry:JCMCC:1998} and a new efficient matching algorithm. This new algorithm turned out to be a crucial part of the discovery of these new matrices. \\

\noindent
A Hadamard matrix of order $428$ has been constructed in \cite{KT:2005:428} and currently the smallest multiple of $4$ for which a Hadamard matrix is not known is $668$. We refer the reader to \cite{Horadam:2007} for more information on Hadamard and skew-Hadamard matrices. \\

\noindent
The updated list of integers $v<500$ for which no Hadamard matrices of order $4v$ are known
consists now of $12$ integers
$$
    167, 179, 223, 283, 311, 347, 359, 419, 443, 479, 487, 491
$$
all of them primes congruent to $3 \,\, (\!\!\!\!\mod 4)$.

\section{The matching algorithm}

\noindent One of the serious difficulties of the unions of orbits approach when used in conjunction
with the Goethals-Seidel array is that one needs to locate four lines in four text files,
each line containing an $n$-tuple of non-negative integers,
such that the element-wise sum of the four $n$-tuples
is an $n$-tuple whose entries are all equal to a pre-defined constant. \\

\noindent The difficulty stems from the fact that if two files contain ten million lines each, then the file
with all possible sum combinations of these two files will contain $10^{14}$ lines, which is impractical even
to store in a file (with $10$ bytes per line, this would take $1$ petabyte).
Therefore the ensuing naive algorithm to solve this problem, i.e. to combine the files two by two in two pairs and look for a potential match, is utterly impractical. In order to circumvent this difficulty, we compress the lines using linear hashing and parallelize the computation, so that each worker processes a small subset of the sums that easily fits in memory. \\

\noindent Let us begin with a precise statement of the problem we were faced with.
Given four text files with $N_1$, $N_2$, $N_3$, $N_4$ lines each,
such that each line contains an $n$-tuple of non-negative integers,
one has to find four lines, one line in each file, identified by their line numbers
$l_1$, $l_2$, $l_3$, $l_4$, such that the element-wise sum of the $n$-tuples in those lines equals an $n$-tuple of the form $( \lambda, \ldots, \lambda )$.
In the case of the SDS $(631;315,330,330,330;674)$ it turns out that $n = 21$, $\lambda = 674$ and that
the four files contained approximately $10$ million lines each. Note that since $330$ appears three times
in this particular SDS, three of the input files are in fact identical, so there are only two different files
in total, among the four files processed by the algorithm. \\

\noindent Let $\Z_+^n$ denote the additive monoid of non-negative integer $n$-tuples.
Given four subsets $A_1,A_2,A_3,A_4\in \Z_+^n$ and a target $s\in \Z_+^n$,
we need to find $a_i\in A_i$, $i=1,2,3,4$, such that $a_1 + a_2 + a_3 + a_4 = s$,
or prove that they do not exist. \\

\noindent A na\"\i ve exhaustive search requires $O(\prod |A_i|)$ time and constant space.
A meet-in-the-middle approach \cite{Horowitz:Sahni:JACM:1974} reduces this to an amortized
$O(|A_1||A_2| + |A_3||A_4|)$ time but takes $O(|A_1||A_2|)$ space,
by first storing all the sums $a_1 + a_2$, where
$a_1\in A_1$, $a_2\in A_2$, in a hash set $H$, and then searching for
$a_3\in A_3$, $a_4\in A_4$, such that $s - a_3 - a_4 \in H$.
The above amortized complexity is due to the fact that if there are not too many collisions then the complexity of
hash table operations is constant time, see \cite{CLRS:2009}.
The hash set with amortized constant time insertions and lookups was implemented as a linear probing \cite{Pagh:Pagh:Ruzic:2009}
hash table with step size $1$, using T. Wang's $64$-bit hash function than can be found on-line
at \url{http://en.wikipedia.org/wiki/Hash_table} for instance.

\noindent For convenience, we took $B_1 = \{\lfloor \frac s2\rfloor - a_3\; |\; a_3 \in A_3\}$ and
$B_2 = \{ \lceil \frac s2 \rceil - a_4\; |\; a_4 \in A_4\}$ and reduced the original
problem to the equivalent one of finding $a_1\in A_1$, $a_2\in A_2$, $b_1\in B_1$,
$b_2\in B_2$, such that $a_1 + a_2 = b_1 + b_2$. In our case, the elements of
$B_1$ and $B_2$ turned out to be non-negative integer tuples. \\

\noindent In order to reduce space requirements and speed up additions, we mapped the $n$-tuples to 64-bit
non-negative integers using a linear hash function $h: \Z_+^n \to \Z/2^{64}\Z$. The linear hash function was
defined as follows. First we subtracted a suitable constant from each entry, so that the result is between
$0$ and $127$. Then, for each line, we partitioned the numbers into groups of $8$, packed each group into
an $8$-byte machine word, and took the linear combination of the resulting machine words with random odd
coefficients (the coefficients were chosen in advance). This hash function has the following property:
for any pair of input lines containing $n$-tuples $a,b$, we have that $h(a+b) = h(a) + h(b)$.
Then we listed the solutions of $h(a_1) + h(a_2) = h(b_1) + h(b_2)$,
and for each of them checked whether they yield a solution to the original problem
(some of them may not because of collisions). \\

\noindent Furthermore, we parallelized the algorithm, reducing the overall search time and
the space requirements on each worker by a factor of $M$, where $M$ is the number
of workers.
First, we represented each of the four input sets (after hashing) as an array of lists,
which at index $i$, $0\le i<M$, stores the list of all elements equal to $i$
modulo $M$.
Using this data structure, worker number $i$, $0\le i < M$, can easily enumerate
and store in the hash set $H$ all sums $h(a_1) + h(a_2)$ equal to $i$ modulo $M$, and then
lookup in $H$ all sums $h(b_1) + h(b_2)$ that are also equal to $i$ modulo $M$. \\

\noindent In practice, $M$ was chosen to be much greater than the available number of workers,
to further reduce memory use.
Each worker pre-loaded the four sets in memory (represented with arrays of lists,
as described above) and sequentially processed several remainders modulo $M$, where
the number of remainders to be processed in a single run was taken large
enough, so that the search time would dominate the pre-loading time, but would not
exceed the maximal allowed duration of a job on the cluster.
Reduced memory requirements enabled us to simultaneously schedule multiple workers on
a single multi-core machine, thus fully utilizing the capacity of the cluster.

\section{Results}

\noindent We now present our results in the form of SDSs, for $v = 213, 251, 631$.
Non-equivalence of SDSs was established by an implementation of the method described in \cite{Djokovic:AnnComb:2011}. \\

\noindent We define the notation for the orbits of the action of
the subgroup that we use to construct the solutions below.
The automorphism group of the additive cyclic group $\Z_v$ will
be identified with the group of invertible elements, $\Z_v^\star$,
of the ring $\Z_v$. Consider a
fixed subgroup $H$ of order $h$ of $\Z_v^\star$. Clearly $h$ must
divide $|\Z_v^\star| = \phi(v)$. Denote by $H \cdot k$ the orbit of $H$ in $\Z_v$ through the point $k$, where $\cdot$ is multiplication $\mod \,\, v$.
We refer to the orbit $H \cdot 0 = \{ 0 \}$ as the
{\it trivial orbit}. The orbit $H \cdot 1$ is just the subgroup
$H$ itself. In general the size of an orbit may be any divisor of
$|H|$ and if $v$ is a prime then every nonzero orbit is just a
coset of $H$ in $\Z_v^\star$ and so will have size $|H|$.

\noindent The notation
\begin{equation}
      X = \bigcup_{j \in J} H \cdot j,
      \quad
      Y = \bigcup_{k \in K} H \cdot k,
      \quad
      Z = \bigcup_{l \in L} H \cdot l,
      \quad
      W = \bigcup_{m \in M} H \cdot m
        \label{Notation:JKLM_SDS}
\end{equation}
will be used below to present all the solutions found, in fact each
solution will be given only via the four index sets $J, K, L, M$.
For suitable choices of the four index sets $J,K,L,M$,
the four sets $X$, $Y$, $Z$, $W$ defined in (\ref{Notation:JKLM_SDS}) form
SDS$(v;x,y,z,w;\lambda)$ with $x = \sum_{j \in J} |H \cdot j|$, $y =
\sum_{k \in K} |H \cdot k|$, $z = \sum_{l \in L} |H \cdot l|$, $w =
\sum_{m \in M} |H \cdot m|$ and $\lambda = x+y+z+w-v$. \\

\noindent All Hadamard matrices in this paper are constructed via the Goethals-Seidel array (see \cite{Seberry:Yamada:1992})
$$
\left[
\begin{array}{cccc}
P_1    &  P_2 R      &  P_3 R    &  P_4 R    \\
-P_2 R &  P_1        &  -P_4^T R &  P_3^T R  \\
-P_3 R &  P_4^T R    &  P_1      &  -P_2^T R \\
-P_4 R &  -P_3^T R   &  P_2^T R  &  P_1      \\
\end{array}
\right]
$$
where $R$ denotes the $v \times v$ matrix with ones in the back-diagonal and zeros everywhere else.
To obtain a Hadamard matrix via the Goethals-Seidel array one has to substitute the $v \times v$ matrices $P_1, P_2, P_3, P_4$ by the $\{ \pm 1 \}$ circulant matrices that arise from the four subsets $X,Y,Z,W$ that make up the SDS$(v;x,y,z,w;\lambda)$. More precisely, let $a_X = (a_0,\ldots,a_{v-1})$ be a sequence defined by $a_i = -1$ if $i \in X$, and $a_i = 1$ if $i \not\in X$; and define the sequences $a_Y, a_Z, a_W$ similarly. Then denote by $[X],[Y],[Z],[W]$
the $\{ \pm 1 \}$ circulant matrices whose first rows are $a_X, a_Y, a_Z, a_W$ respectively.
If $X,Y,Z,W$ form a SDS$(v;x,y,z,w;\lambda)$, then
$$
        [X][X]^T + [Y][Y]^T + [Z][Z]^T + [W][W]^T = 4v I_v.
$$

The parameter sets $(v;k_1,k_2,k_3,k_4;\lambda)$ that we need can be constructed as follows. Assuming that $v$ is odd, we find all representations of $4v$ as a sum of four odd squares
\begin{equation} \label{zbir4kv}
4v=n_1^2+n_2^2+n_3^2+n_4^2
\end{equation}
with $n_1\ge n_2\ge n_3\ge n_4>0$ and $n_1<v/2$. We set $k_i=(v-n_i)/2$. Then one can verify that $\sum k_i(k_i-1)=\lambda(v-1)$ where $\lambda=\sum k_i -v$. This assertion remains valid if instead of $k_i=(v-n_i)/2$ we set $k_i=(v+n_i)/2$ for some indexes $i$. This is useful because we need each base block $X_i$ to be a union of orbits of $H$. For instance if $v$ is a prime number, and $h$ is the order of the subgroup $H$, then each nontrivial orbit has size $h$ and so each $k_i$ must be divisible by $h$. However, $h$ may divide $(v+n_i)/2$ but not $(v-n_i)/2$. We say that the decomposition (\ref{zbir4kv}) of $4v$ into sum of four odd squares is associated with the parameter set
$(v;k_1,k_2,k_3,k_4;\lambda)$. \\

Suppose now that we have already selected a suitable parameter set $(v;k_1,k_2,k_3,k_4;\lambda)$ and a suitable subgroup $H$ such that each $k_i$ is a sum of the cardinalities of certain
$H$-orbits. Since we build the base blocks $X_i$ from the nontrivial orbits of $H$, the most important quantity is the number $\nu$ of these orbits. If $v$ is a prime number, then $\nu$ is just the index of $H$ in $\Z_v^*$. In the cases that we succeeded to construct SDSs $\nu$ did not exceed 50.
For each $i=1,2,3,4$ we construct two files: the file $F_i$ such that each line of the file contains the list of the labels of orbits whose union $X_i$ has cardinality $k_i$ and the file $F'_i$ listing the multiplicities of the differences $x-y\in\Z_v$ for all ordered pairs $(x,y)$, with $x,y\in X_i$ and $x\ne y$. Apart from the cases where $\nu$ is small, it is impossible to carry out an exhaustive search for all possible candidates for the base block $X_i$ and so we start the search at a randomly chosen place and run our program to collect the desired number of candidates. After constructing these 8 files, we have to find a match in the four files $F'_i$. By a ``match'' we mean that the sum of the multiplicities from the four suitably selected lines from the $F'_i$, one line from each file, is constant (necessarily equal to $\lambda$). Once the match is found we construct the base blocks $X_i$ by using the files $F_i$ and the line numbers provided by the match.

\subsection{$v = 213$}
Consider the subgroup $H = \{ 1,37,91,103,172,187,190 \}$ of order $7$, of $\Z_{213}^\star$.
Note that $\Z_{213}^\star$ is of order $\phi(213) = 140$.
There are $32$ nontrivial orbits of the action of $H$ on $\Z_{213}$.
We give an SDS with parameters $(213;106,106,105,92;196)$,
via the index sets $J,K,L,M$ (with respective cardinalities $16,15,14,16$) to be used in
(\ref{Notation:JKLM_SDS}), which gives rise to skew-Hadamard matrices
of order $4 \cdot 213 = 852$. The associated decomposition into  the sum of $4$ squares is:
$$
\begin{array}{ccl}
    4 \cdot 213 = 852 & = & 1^2 + 1^2 + 3^2 + 29^2 = \\
                       & = & (213-2 \cdot 106)^2 + (213-2 \cdot 106)^2 + (213-2 \cdot 105)^2 +(213-2 \cdot 92)^2.
\end{array}
$$
$$
\begin{array}{l}
\JKLM{4, 5, 7, 10, 11, 15, 17, 19, 20, 30, 34, 38, 39, 42, 43, 142}{2, 7, 12, 14, 17, 22, 28, 34, 39, 42, 43, 44, 69, 84, 86}
{1, 4, 11, 15, 17, 20, 21, 28, 30, 34, 42, 44, 69, 142}{2, 4, 10, 11, 12, 15, 21, 22, 23, 28, 30, 34, 44, 69, 71, 86} \\
\end{array}
$$

\subsection{$v = 251$}
Consider the subgroup $H = \{ 1,20,113,149,219 \}$ of order $5$, of $\Z_{251}^\star$.
Note that $\Z_{251}^\star$ is of order $\phi(251) = 250$.
There are $50$ nontrivial orbits of the action of $H$ on $\Z_{251}$, all of size $5$.
We give two SDSs with parameters $(251;125,120,115,115;224)$,
via the index sets $J,K,L,M$ (with respective cardinalities $25,24,23,23$) to be used in
(\ref{Notation:JKLM_SDS}), which give rise to Hadamard matrices
of order $4 \cdot 251 = 1004$. The associated decomposition
into the sum of $4$ squares is:
$$
\begin{array}{ccl}
    4 \cdot 251 = 1004 & = & 1^2 + 11^2 + 21^2 + 21^2 = \\
                       & = & (251-2 \cdot 125)^2 + (251-2 \cdot 120)^2 +
                             (251-2 \cdot 115)^2 + (251-2 \cdot 115)^2.
\end{array}
$$
$$
\begin{array}{l}
\JKLM
{2, 4, 5, 6, 7, 111, 9, 37, 10, 11, 12, 14, 173, 15, 24, 74, 19, 48, 33, 73, 53, 57, 43, 72, 68}
{2, 55, 6, 75, 16, 11, 173, 24, 17, 42, 18, 19, 106, 21, 30, 48, 34, 73, 35, 53, 41, 57, 43, 97}
{1, 2, 4, 50, 6, 75, 9, 37, 10, 11, 173, 15, 17, 18, 21, 48, 36, 34, 119, 45, 72, 68, 97}
{32, 2, 25, 3, 5, 8, 7, 16, 173, 15, 74, 19, 106, 21, 95, 30, 73, 53, 41, 57, 45, 72, 68} \\
\end{array}
$$
$$
\begin{array}{l}
\JKLM
{32, 3, 55, 50, 5, 6, 7, 9, 37, 11, 12, 14, 15, 24, 18, 19, 82, 30, 36, 34, 73, 35, 41, 57, 72}
{1, 32, 25, 55, 5, 8, 75, 9, 37, 11, 14, 15, 24, 17, 42, 18, 19, 95, 48, 41, 57, 43, 45, 97}
{1, 2, 25, 3, 50, 8, 75, 9, 11, 12, 173, 24, 17, 28, 95, 30, 33, 36, 34, 73, 45, 72, 97}
{32, 25, 4, 5, 6, 75, 7, 16, 11, 12, 24, 42, 106, 28, 95, 48, 73, 35, 53, 41, 57, 45, 72} \\
\end{array}
$$

\subsection{$v = 631$}

Consider the subgroup $H = \{ 1, 8, 43, 64, 79, 188, 228, 242, 279, 310, 339, 344, 512, 562, 587 \}$
of order $15$, of $\Z_{631}^\star$. There are $42$ nontrivial orbits of the action of $H$ on $\Z_{631}$
all of size $15$.
We give four nonequivalent SDSs with parameters $(631;315,330,330,330;674)$,
via their index sets $J,K,L,M$ (with respective cardinalities $21,22,22,22$) to be used in
(\ref{Notation:JKLM_SDS}), which give rise to Hadamard matrices
of order $4 \cdot 631 = 2524$. In addition, the first two SDSs give rise to skew-Hadamard matrices of
order $4 \cdot 631 = 2524$. The associated decomposition into
the sum of $4$ squares is:
$$
\begin{array}{ccl}
    4 \cdot 631 = 2524 & = & 1^2 + 29^2 + 29^2 + 29^2 = \\
                       & = & (631-2 \cdot 315)^2 + (631-2 \cdot 330)^2 + (631-2 \cdot 330)^2 +(631-2 \cdot 330)^2.
\end{array}
$$

$$
\JKLM{ 1, 2, 3, 4, 6, 7, 12, 13, 14, 17, 19, 21, 26, 27, 31, 38, 42, 52, 62, 76, 124 }{ 1, 2, 3, 4, 5, 6, 7, 9, 12, 17, 23, 26, 27, 31, 33, 38, 42, 44, 52, 76, 78, 126 }{ 1, 2, 3, 4, 6, 7, 9, 11, 12, 13, 14, 17, 18, 19, 21, 29, 35, 46, 52, 62, 66, 76 }{ 1, 2, 3, 4, 5, 6, 12, 13, 14, 18, 19, 21, 22, 26, 27, 38, 39, 42, 63, 67, 92, 124 }
$$
$$
\JKLM{ 11, 13, 19, 22, 26, 29, 31, 33, 38, 39, 44, 52, 62, 65, 66, 67, 76, 78, 117, 124, 187 }{ 1, 3, 4, 5, 9, 11, 14, 17, 18, 22, 23, 26, 29, 33, 38, 39, 42, 46, 62, 65, 67, 117 }{ 1, 2, 5, 6, 9, 17, 18, 21, 22, 27, 33, 39, 44, 46, 52, 66, 76, 78, 92, 117, 124, 187 }{ 2, 4, 5, 7, 9, 11, 12, 13, 18, 19, 21, 23, 29, 31, 42, 44, 65, 66, 67, 78, 92, 187}
$$
$$
\JKLM{ 1, 3, 4, 7, 12, 13, 17, 18, 19, 27, 29, 31, 33, 35, 42, 46, 62, 67, 92, 124, 187 }{ 5, 6, 7, 9, 12, 17, 18, 21, 22, 26, 29, 33, 35, 38, 42, 44, 62, 63, 67, 76, 124, 126 }{ 4, 6, 7, 9, 11, 13, 14, 17, 18, 21, 22, 23, 27, 29, 35, 39, 46, 67, 78, 124, 126, 187 }{ 2, 5, 6, 7, 11, 13, 14, 18, 19, 21, 27, 29, 33, 35, 39, 44, 52, 62, 63, 65, 117, 187 }
$$
$$
\JKLM{ 1, 2, 4, 5, 6, 9, 11, 12, 13, 17, 22, 27, 29, 42, 44, 65, 78, 92, 117, 126, 187 }{ 1, 7, 9, 13, 17, 19, 22, 23, 27, 29, 31, 33, 35, 38, 44, 46, 65, 66, 76, 78, 126, 187 }{ 2, 3, 4, 5, 13, 17, 21, 22, 29, 35, 38, 39, 52, 62, 63, 65, 67, 76, 92, 117, 124, 126 }{ 2, 3, 4, 5, 7, 9, 11, 12, 17, 18, 21, 22, 29, 33, 42, 44, 46, 52, 65, 66, 92, 187}
$$

\noindent We now report matching algorithm timings\footnote{See \url{//computing.llnl.gov/tutorials/parallel_comp/} for the definition of the term ``task''.} for each of the three cases $v = 213, 251, 631$.
\begin{itemize}
\item $v=213$: Each input file contained $10$ million lines. The matching algorithm launched $1$ single-core task which ran for less than one minute before a match was found.

\item $v = 251$: Each input file contained $10$ million lines. The matching algorithm launched $140$ single-core tasks which ran in parallel for $40$ minutes, until one of them found a match.

\item $v = 631$: Each input file contained $10$ million lines. The matching
algorithm launched $325$ single-core tasks, which ran in parallel for $8$
minutes, until one of them found a match.

\end{itemize}
Since the tasks ran independently of each other and
the startup time for each task was negligible, the speedup was linear.
The computations were performed on a RQCHP supercomputer with 1588 computing nodes SGI C2112-4G3
with the following characteristics: 2 AMD processors 12 cores 6172 2.1GHz
32 GB of RAM memory and an 1TB Hard Disk.

\section{Acknowledgements}
The authors thank the referees for their pertinent comments and suggestions.
The first and the third author wish to acknowledge generous support by NSERC.
This work was made possible by the facilities of the Shared Hierarchical
Academic Research Computing Network (SHARCNET) and Compute Ontario, as
well as the facilities of the R\'eseau Qu\'eb\'ecois de Calcul de Haute Performance (RQCHP)
and Calcul Qu\'ebec.

\end{document}